\documentstyle[12pt]{article}
\begin{document}

\thispagestyle{empty}
\vspace*{0.4cm}

\hfill {\begin{tabular}{c} q-alg/9512017 \\
December 1995 \end{tabular} }

\vspace*{1.6cm}  \begin{center}
{\Large \bf QUANTUM GROUP \\
\vspace*{2mm}
COVARIANT SYSTEMS }  \end{center}
\vspace{1cm}

\begin{center}  {\bf M. Chaichian} $^*$  \end{center}
\vspace*{0.1cm}

\begin{center} \begin{tabular}{c}
\it High Energy Physics Laboratory, Department of Physics, \\
\it Research Institute for High Energy Physics \\
\it P.O. Box 9 (Siltavuorenpenger 20C), University of Helsinki  \\
\it Helsinki SF-00014, Finland
\end{tabular} \end{center}
 \vspace*{0.4cm}

\begin{center}  {\bf P.P. Kulish} $^{**}$  \end{center}
\vspace*{0.1cm}

\begin{center} \begin{tabular}{c}
\it St.Petersburg Branch of Steklov Mathematical Institute \\
\it Fontanka 27, St.Petersburg, 191011, Russia  \\
\end{tabular} \end{center}

 \vspace*{1.5cm}

\begin{center} \bf{ Abstract  } \end{center}
The meaning of quantum group transformation properties is discussed
in some detail by comparing the (co)actions of the quantum group with
those  of the corresponding Lie group, both of which have the same
algebraic (matrix) form of the transformation. Various algebras are considered
which are covariant with respect to the quantum (super) groups
$SU_q(2),\; SU_q(1, 1),\; SU_q(1|1),\; SU_q(n), \\
SU_q(m|n),\; OSp_q(1|2)$ as well as deformed Minkowski space-time algebras.

\vfill

\begin{flushleft}   \rule{6 cm}{0.05 mm} \\
$^{*}$ {\footnotesize e-mail:$\;$ chaichian@phcu.helsinki.fi } \\
Š
$^{**}$ {\footnotesize e-mail:$\;$ kulish@pdmi.ras.ru } \\
\end{flushleft}

\newpage

\begin{center}
{\LARGE {\bf QUANTUM GROUP COVARIANT SYSTEMS}}
\end{center}

\vspace{0.4cm}

\begin{center}
{\Large {M. Chaichian and P.P. Kulish}}
\end{center}

\section{Introduction}
The transformation properties of physical systems related to the Lie
groups are of great importance for the understanding of Nature.
As a result, applications of the Lie group theory  take place
in quite different branches of
physics and the corresponding formalism is very well developed. The quantum
groups and quantum algebras extracted from the quantum inverse scattering
method (QISM) happen to be quite similar or even richer mathematical objects
as compared to Lie groups and Lie algebras.

In this paper we shall point out some peculiarities of the quantum group
interpretations when the formal transformations of the physical
quantities coincide with the usual
ones while the coefficients (elements) of these transformations
are now non-commutative quantities belonging to a quantum group (QG)
or a quantum algebra.

These objects (QG and q-algebras) are described using the language of
Hopf algebras. In the general situation of Lie group theory one has the Lie
algebra $Lie (G)$ (or better to say the corresponding universal enveloping
algebra) with non-commutative multiplication and symmetric coproduct
$\Delta$, and the commutative algebra of functions $F$ on the Lie group
manifold $G$ with non-symmetric coproduct
$\Delta : F \rightarrow  F \otimes F  $. After a q-deformation (or
"quantization") the corresponding objects $Lie_{q} (G)$ and $F_q$ start to be
much more similar in between as the Hopf algebras with non-commutative
multiplications and non-symmetric coproducts in both cases. Hence, it looks
natural to have the same physical interpretation
for transformations including both of them.

Putting aside the complicated integrable models solved by the quantum inverse
scattering method and the quantum conformal field theory we consider systems
with finite degrees of freedom such as a set of q-oscillators
${\cal A}$ covariant under the coaction $\varphi$ of the
quantum (super-) group $F_q$.
When the coaction $\varphi: {\cal A} \rightarrow  F_q  \otimes {\cal A}$
does not preserve the physical observables such as
Hamiltonian, momenta, etc, the standard problem of the tensor product
decomposition of ${\cal H}_F \otimes {\cal H}_{\cal A}$ emerges, where
${\cal H}_F$ and ${\cal H}_{\cal A}$ are the state spaces of the corresponding
algebras. If instead of $F_q$ we have the quantum algebra $Lie_{q} (G)$
then its representations are almost identical to the undeformed algebra (let
us omit plenty of technicalities related with the case when $q$ is a
root of unity). However,
the deformed algebra of functions $F_q$ is a new algebraic object with
its own representation theory and the corresponding decomposition problem
of ${\cal H}_F \otimes {\cal H}_{\cal A}$ or
${\cal H}_F \otimes {\cal H}_F$  is new as well.

\section{Covariant systems}

2.1. ${\cal A}_{q}$ as $su_{q}(2)$-covariant algebra.

Let us start with a simple covariant system.
The q-oscillator algebra ${\cal A}_{q}$ has three generators
with commutation relations (for $q$ real one has
$ A = q^{N/2} a = q^{N} \alpha $)

\begin{equation}
a a^{ \dagger} - q  a^{ \dagger} a = q^{-N}, \quad  [N , a ] = - a \; ,
\quad  [N , a^{ \dagger} ] =  a^{ \dagger} \; ;
\end{equation}

\begin{equation}
A A^{ \dagger} - q^2 A^{ \dagger} A  = 1,  \quad [N, A ] = - A \;, \quad
[N, A^{ \dagger} ] =  A^{ \dagger} \; ;
\end{equation}

\begin{equation}
[\alpha , \alpha^{ \dagger}] = q^{-2N},  \quad [N, \alpha ] = - \alpha , \quad
[N, \alpha^{ \dagger} ] = \alpha^{ \dagger} \;.
\end{equation}

\noindent
The third set can be obtained  from
the quantum algebra $su_{q}(2)$ (with generators $X_{+}, X_{-}, J$ and
well-known commutation relations and coproduct $\Delta$)
by a contraction procedure with fixed $q$ \cite{8} ($\lambda = q- q^{-1} $)

$$
\alpha = lim \lambda^{1/2}X_{+}/q^s , \quad s \rightarrow \infty \;,
\quad N = s - J \;.
$$

\noindent
{}From this contraction one could find that although
there is no natural coproduct for ${\cal A}_{q}$,
the formulas for the $su_{q}(2)$ coproduct that survive under
the contraction procedure could be
interpreted as covariance of the algebra ${\cal A}_{q}$
with respect to the quantum algebra $su_{q}(2)$.
The corresponding map is
$\psi : {\cal A}_{q} \rightarrow {\cal A}_{q} \otimes su_{q} (2)$, such that
$$
\psi (N) = N - J,
$$
\begin{equation}
\psi (\alpha ) = \alpha q^{-J} + \sqrt{\lambda} q^{-N} X_{+},
\end{equation}
$$
\psi (\alpha^{\dagger})= \alpha^{\dagger} q^{-J} + \sqrt{\lambda}
q^{-N} X_{-}.
$$

\noindent
It is easy to check the following consistency properties for this coaction :
$(\psi \otimes id) \circ \psi = (id \otimes \Delta ) \circ \psi$ and
$(id \otimes \epsilon ) \circ \psi= id$ as well as that
$\psi$ preserves the defining
relations of ${\cal A}_{q}$. However, the central element $z$
of the algebra ${\cal A}_{q}$

\begin{equation}
z= \alpha^{\dagger} \alpha - [N; q^{-2}]\;, \quad
[x; q] = (1 - q^x)/(1-q)\;,
\end{equation}

\noindent
is not invariant under this coaction: $\psi (z) \neq z$.

If we choose the Hamiltonian of the $q$-oscillator as $H=\alpha^{\dagger}
\alpha$ and restrict ourselves to the irreducible representation
${\cal H}_{F}$ of ${\cal A}_{q}$
with the vacuum state: $\alpha |0>=0, N|0>=0$ (for $q \in (0, 1)$ there are
other irreps \cite{5,8}) then $z=0$ and

\begin{equation}
H= \alpha^{\dagger} \alpha = [N; q^{-2}] = (1 - q^{-2N})/(1-q^{-2}).
\end{equation}

\noindent
The spectrum of $H$ and its eigenstates are obvious

$$
spec\; H= \{ [n; q^{-2}]\;,\; n=0, 1, 2, \ldots \},
$$

$$
|n>= ([n; q^{-2}]!)^{-1/2} (\alpha^{\dagger})^{n} |0>.
$$
\noindent
After the coaction the changed Hamiltonian describes an "interaction" of the
$q$-oscillator with the $q$-spin

\begin{equation}
H_{I} = \alpha^{\dagger} \alpha q^{-2J} + \lambda q^{-2N} X_{-} X_{+}
+ \sqrt{\lambda} q^{-N-J} (X_{-}\alpha /q + q \alpha^{\ast} X_{+})\;.
\end{equation}

\noindent
It acts in the space of physical states
${\cal H}_{ph} = {\cal H}_{F} \otimes V_{j}$,
where $V_{j}$ is an irreducible
finite dimensional representation of the $su_{q}(2)$ of dimension $2j+1$. This
space is decomposed into the direct sum of $2j + 1$ irreducible
representations of ${\cal A}_{q}$: ${\cal H}_{ph} = \sum_k {\cal H}_{F}^{(k)}$
with corresponding vacuum states $|0>_{k} \, , \, k=0, 1,\ldots , 2j$

\begin{equation}
|0>_{k} =\sum^{k}_{m=0} |m> \otimes |j-k+m; j>c(m,k) \; ,
\end{equation}

\noindent
where
$$
|m> \in {\cal H}_{F} \; , \; |l ; j > \in V_{j} \; ,
\; J|l ; j> = l |l ; j > ,
$$

$$
X_{+}|m \; ;\; j > = N_+(m,j)|m+1 \; ;\; j > \; , \;
N_+^2(m,j) = [j-m]_q [j+m+1]_q \;,
$$
$$
c(m,k)= (- \sqrt{\lambda} q^{j-k+1} )^m ([m;q^{-2}]!)^{-1/2}
\prod_{l=0}^{m-1} N_+(j-k+l,\;j)\;.
$$

\noindent
The spectrum of $H_{I}$ coincides
with that of $H$ up to the multiplicative factor $q^{2(j-k)} $
in each subspace $ {\cal H}_{F}^{(k)}$, but it has the multiplicity $2j + 1$
$$
spec\; H = \{ q^{2(j-k)} [n; q^{-2}]\;,\; n=0, 1, 2, \ldots \} .
$$

The central element
$\psi (z)$ has $2j+1$ eigenvalues $-[k-j ; q^{-2}]$.

It is interesting to point out that this coaction has no classical
(non-deformed) counterpart in the limit $q \rightarrow 1$ in the quantum
theory. Such a limit exists in the Poisson-Lie theory (see e.g. \cite{B}).
The connection of the $q$-oscillator algebra ${\cal A}_{q}$ with
$su_q(2)$ through the contraction procedure gives rise also
to a more complicated coaction of the quantum group
$SU_q(2)$ on ${\cal A}_{q}$ (see Subsec. 2.5 and 2.6).

2.2. ${\cal A}_{q}$ as $SU_{q} (1, 1)$-covariant algebra.

Let us consider the second set (2) of the $q$-oscillator algebra
${\cal A}_{q}$ generators with relation
$$
A A^{\dagger}=q A^{\dagger} A + 1 .
$$

\noindent
redenoting $q^{2}$ by $q$ and putting $q \in (0,1)$.
This relation reminds us of the quantum plane with $xy=qyx$ a central
extension. Using the two component column $X^{t}= (A, A^{\dagger})$
it can be rewritten in the $R$-matrix form \cite{7,V}

\begin{equation}
\hat{R} X \otimes X = q X \otimes X - q^{-1} J,
\end{equation}

\noindent
where $\hat{R} = {\cal P}R$ is the $R$-matrix of $su_{q}(2)$ and
$J$ is the four
component column $J^{t}= (0, 1, -q, 0)$ obviously related to the well-known
$q$-metric $2 \times 2$ matrix $\epsilon_{q}$. This relation is preserved
under the transformation $\psi : X \rightarrow \psi (X)= TX$, with $T$
being the $2 \times 2$ matrix of the quantum group $SU_{q} (1, 1)$ generators
$$
T=
\left(
\begin{array}{ll}
a & b\\
b^{\ast} & a^{\ast}
\end{array}
\right)
$$

\noindent
which satisfies the FRT-relation $\hat R T\otimes T = T\otimes T \hat R $ [1].
The invariance of the inhomogeneous term is just another form of the
$q$-metric relation

\begin{equation}
T \epsilon_{q} T^{t} = \epsilon_{q} det_{q} T, \quad T \otimes T J =
T_1 T_2 J = det_{q}T J \;,
\end{equation}

\noindent
provided that the quantum determinant of $T$ is $1$:
$det_{q} T = a a^{\ast} - q b b^{\ast} = a^{\ast} a - b^{\ast}b/q=1$
(a defining condition for $SU_{q}(1, 1)$). (One can consider central
extension of the real quantum plane as well with $|q|=1$ covariant with
respect to $SL_{q}(2, {\bf R})$. Then the reality condition will fix the
phase of the constant term.) The map $\psi : {\cal A}_{q} \rightarrow
{\cal A}_{q} \otimes SU_{q} (1, 1)$ or in terms of the generators

\begin{equation}
\psi (A) = aA + b A^{\dagger} \; , \; \psi (A^{\dagger})= a^{\ast}
A^{\dagger} + b^{\ast} A
\end{equation}

\noindent
satisfies all properties of a coaction. Its form is reminiscent of the famous
Bogoliubov transformation. However, now the "coefficients" are
non-commuting. The $q$-oscillator Hamiltonian acting in the same
space ${\cal H}_{F}$ as (6)

\begin{equation}
H= A^{\dagger}A= [N; q]= (q^{N}-1)/(q-1)
\end{equation}

\noindent
(it differs from the previous one by $q^{N}$ factor and renotation of
$q^{2}$) is also not invariant under the coaction $\psi$

\begin{equation}
\psi (H)= A^{\dagger} A + b^{\ast} b + [2]_{q} A^{\dagger} A b^{\ast} b +
b^{\ast} a A^{2} + a^{\ast} b (A^{\dagger})^{2}.
\end{equation}

\noindent
The commutation relations of the $SU_{q}(1, 1)$ generators (as well as those
of $SU_{q}(2)$)
$$
ab = qba \quad , \quad a b^{\ast} = q b^{\ast} a, \quad
b b^{\ast} = b^{\ast} b \;,
$$
$$
[a, a^{\ast}] = \lambda b^{\ast} b, \quad \quad \lambda = (q-1/q)\;,
$$

\noindent
themselves remind us of the $q$-oscillator algebra $(a \sim \alpha^{\dagger},
a^{\ast} \sim \alpha , b \sim b^{\ast} \sim q^{-N}$) with the additional
condition due to $det_{q} T=1$
$$
a a^{\ast} = 1 + q b^{\ast} b, \quad a^{\ast} a = 1 + b^{\ast} b/q .
$$

\noindent
The deformation parameter $q$ being less than 1 forces us to consider the
irreducible representation of $SU_q(1, 1)$ in the Hilbert
space $l_{2} ({\bf Z})$ with basis
$|m>, m =\ldots , -2, -1, 0, 1, 2. \ldots $ which consists of the
eigenstates of commuting operators $b, b^{\ast}$ for which  $a$
acts as a creation (shift) operator
$$
b|m>= e^{i \phi} q^{-m} |m>, \; a|m> = c_{m} |m+1>,
a^{\ast} |m>= c_{m-1} |m-1>,
$$
$$
c_{m}^{2} = 1 + q^{-2m-1}.
$$

\noindent
Hence the transformed Hamiltonian is defined in the space ${\cal H}_{F} \otimes
l_{2} ({\bf Z})$. It has the same spectrum $\{[n; q], \; n=0, 1, \ldots \}$
with infinite multiplicity. The corresponding vacuum states are

\begin{equation}
|0>_{(\psi)}^{(k)} = \sum_{j=0}^{\infty} (-1)^{j} x_{j} |2j> \otimes
|k-j>,
\end{equation}

$$
x_{j}=([2j-1;q]!!/[2j;q]!!)^{1/2} \Pi_{i=1}^{j} q^{i-k}/c_{k-i} \; ,
$$

\noindent
where the second vector in the tensor product belongs to the $SU_{q} (1, 1)$
irrep space $l_{2} ({\bf Z})$.

As in the case of the representation theory the {\it invariant subspaces} of
the QG $F$ (coaction) corepresentation $V$ can be defined as
$W \subset V$ such that $\phi : W \rightarrow F \otimes W$ .
%In our initial example of the quantum plane the linear space $z_1x + z_2y$ is
%irreducible, while the whole covariant algebra $C_q^2$ is decomposed into
%irreducible subspaces
%$C_q^2= \oplus_s V_s , s= 0, 1/2, 1, ...; dimV_s=2s+1 $.
%The subspaces $V_s$ consist of linear combinations of
%$2s$-powers monomials in
%$x$ and $y$ and generate covariant subalgebras ${\cal A}_s$ such that
%$ {\cal A}_s \subset C_q^2 $, but not $ {\cal A}_{s+1} \subset {\cal A}_s  $.
The {\it invariant elements} of the $F$-corepresentation $V$
do not change at all:
$\phi (v)= v$ (or better to write $1_F \otimes v$ for one has the possibility
in the corresponding representation of the dual Hopf algebra $(F)^*$ to
contract a dual element $X \in (F)^* $ with $1_F$ to get a number $X(1_F)$ ).

The extensions of the previous examples to higher rank quantum groups
give rise to covariant algebras corresponding to different quantum
homogeneous spaces \cite{VS}, systems of (super)
$q$-oscillators \cite{5,4,6,13}
and examples of non-commutative geometry.

2.3. The covariant super-q-oscillator algebra $s$-${\cal A}_q$ [4]
refers to the quantum super-group
$S U_{q} (1 | 1)$, with the $T$-matrix of the generators
$$
T =
\left(
\begin{array}{cc}
a & \beta \\
\gamma & d
\end{array}
\right)
$$

\noindent
and the commutation relations
$$
\begin{array}{cc}
a \beta = q \beta a, & a \gamma = q \gamma a, \\
d \beta = q \beta d, & d \gamma = q \gamma d,
\end{array}
$$
$$
\beta \gamma = - \gamma \beta, \quad \beta^2 = \gamma^2 = 0, \quad
[a,d] = (q - 1 / q)\gamma \beta .
$$

\noindent
Fixing the central element (super-determinant)
$sdet_{q} T = (a - \beta \frac{1}{d} \gamma) / d = 1 $\\
one gets a simple relation between the even generators $d$ and $a$
$$
d=a-\beta \gamma /qd = a-\beta \gamma /qa
$$

\noindent
due to the nilpotency of the odd generators. The involution
(*-operation) is introduced as follows:
$$
d= 1/a^*,  \quad \gamma = d \beta^* d
= a \beta^* a \; .
$$

\noindent
This involution leads to $TT^{\dagger} = 1$ and it is consistent with
the $Z_2$-grading. One has for
the generators $a,\; a^*,\; \beta ,\; \beta^*$
$$
a^* \beta = \beta a^*/q \;, \quad [a, a^*] = (1-q^2) \beta^* \beta \;, \quad
\beta \beta^* = -q^2 \beta^* \beta \;,
$$
$$
a a^* = 1+ \beta \beta^* \; ,   \quad a^* a = 1- \beta^* \beta \;.
$$

\noindent
Introducing $a' = a (1+ \beta^* \beta/2)
= (1-\beta \beta^*/2) a $, one gets $a' a'^* = a'^* a' = 1$ and the
factorization of the $T$-matrix ($a' = 1/ a'^*$ )
$$
T =
\left(
\begin{array}{cc}
a' & 0 \\
0 & 1/a'^*
\end{array}
\right)
\left(
\begin{array}{cc}
(1- \beta^* \beta/2)  & \beta/qa' \\
\beta^* /a'^* &(1 + \beta^* \beta/2)
\end{array}
\right)
$$

\noindent
with unit super-determinant. One concludes that the $q$-deformation
(quantization) of the $SU(1|1)$ super-group is realized by the unitary
scaling operator $\Lambda ,\; \Lambda^* = \Lambda^{-1} $
acting on the Grassmann variablies $\eta$ and  $\eta^*\;, $
which are not quantized $\beta = \Lambda \eta \;,$
$$
\Lambda \beta = q\beta \Lambda \; ,  \quad
\Lambda \beta^* = q\beta^* \Lambda \; , \quad \Lambda 1 = 1
$$

\noindent
and $ a' = exp(i \varphi )\Lambda $. Hence, like in the non-deformed case
the representations of the $SU_q(1|1) $ are parametrized by the phase and
the Grassmann variable.

The corresponding $SU_q(1|1) $ covariant system of
the super-q-oscillator $s$-${\cal A}_q$
has four generators $A, \; A^{\dagger}, \;B, \; B^{\dagger}$
with relations [4]
$$
A A^{\dagger} - q^2 A^{\dagger} A = 1 \; ,  \quad
B B^{\dagger} + B^{\dagger} B = 1 + (q^2 - 1) A^{\dagger} A \;,
$$
$$
A B = q B A \; , \quad A B^{\dagger} = q B^{\dagger} A \; , \quad
B^2 = B^{\dagger 2} = 0 \;.
$$

Using the $R$-matrix formalism similar to (9) it is not difficult to show that
these relations as well as
the Hamiltonian $H = A^{\dagger} A + B^{\dagger} B$ are invariant w.r.t. the
coaction
$$
\left(
\begin{array}{c}
A \\
B
\end{array}
\right)
\rightarrow
\left(
\begin{array}{cc}
a & \beta \\
\gamma & d
\end{array}
\right)
\left(
\begin{array}{c}
A \\
B
\end{array}
\right)
{}.
$$

However, one can consider Hamiltonians which are not invariant w.r.t. the
QG transformation. The latter one will extend the initial system after the
transformation. In particular, one can consider different versions of the
$q$-deformed  $N=2$-SUSY-algebra taken as the super-charges
$$Q = A^{\dagger} B\;,\; Q^{\dagger}= B^{\dagger} A \quad or \quad
Q = A^{\dagger} f\;,\; Q^{\dagger}= A f^{\dagger}\;, $$
where $f$ and $f^{\dagger} $ are free fermions $f = q^{-N} B$
commuting with $A\;,\;A^{\dagger} $ [4]. In all the cases the coaction
does not extend the space of states for the representation
theory of $SU_q(1|1)$
is rather poor and only additional Grassmann parameters appear
after the coaction.
%These constructions can be realized in terms of (commuting and non-commuting)
%q-path integrals as well.

2.4 $SU_q(n)$- and $SU_q(m|n)$-covariant (super) algebras.

Let us introduce $2n$ generating elements of the $SU_q(n)$-covariant
oscillator
algebra ${\cal A}_q(n)$ \cite{5}, written as $n$-component
column  and row vectors
$$
A^t=(A_1,..., A_n), \quad A^{\dagger} =(A_1^{\dagger},... , A_n^{\dagger}).
$$

\noindent
Their commutation relations in the $R$-matrix form (a spectral parameter
independent Zamolodchikov - Faddeev algebra) \cite{7} ($\hat{R}_{\cal P} =
{\cal P} \hat{R}^t {\cal P}$)
$$
\hat{R} A \otimes A = q A \otimes A, \quad
A^{\dagger} \otimes A^{\dagger} \hat{R}_{\cal P} =
q A^{\dagger} \otimes A^{\dagger},
$$
$$
A \otimes A^{\dagger} = q A^{\dagger}_2 \hat{R} A_2 + I
$$

\noindent
demonstrate easily that, due to the FRT-relation $[\hat{R}, T_1T_2]=0$
and $T T^{\dagger} = I$ the coaction
$\phi (A)= T A,\; \phi (A^{\dagger}) = A^{\dagger}T^{\dagger}$
satisfies all the requirements.

One can rewrite these relations in the form (14) using
the $2n$-component vector
$X=(A_1,..., A_n; A^{\dagger}_1,..., A^{\dagger}_n)$ and the corresponding
$2n \times 2n$ matrix $R$, which happens to be the $R$-matrix of the quantum
group $Sp_q(2n)$. Then the inhomogeneous term ($2n$-component vector $J$) is
expressed using the invariant matrix $C$ of $Sp_q(2n)$: $T^tCT=C$ \cite{1,V}.

The invariant Hamiltonian w.r.t. the $SU_q(n)$-coaction is
$$
H = A^{\dagger}_1 A_1 + A^{\dagger}_2 A_2 + \ldots + A^{\dagger}_n A_n \;,
$$
which in the Fock space can be written in terms of the mode number operators
$N_k\;, k=1,\;2,\; \ldots $
$$
H= (q^{2(N_1 +N_2+ \ldots +N_n)} - 1)/(q^2 - 1)\;.
$$
It was already pointed out that the $SU_q(1, 1)$ quantum group is related to
the $q$-oscillator algebra.
The same is true for the $SU_q(2)$: its defining relations
coincide with (3) up to notations and some factors. A realization of the
$SU_q(n)$ requires $n(n-1)/2$ $q$-oscillators and $n-1$ phase
factors \cite{12}. Hence, transforming the algebra
${\cal A}_q(n) \rightarrow SU_q(n) \otimes {\cal A}_q(n) $
one jumps from the $n$ degrees of freedom to the $n(n+1)/2$ ones.

The situation for the quantum super-group $SU_q(m|n)$ is similar. The
covariant super-algebra $s$-${\cal A}_q(m|n)$ has $m$ boson
and $n$ fermion mutually non-commuting $q$-oscillators \cite{4}
or vice versa. The realization of the $SU_q(m|n)$ requires
$ m(m-1)/2 + n(n-1)/2$ $q$-oscillators and $m \times n$ Grassmann parameters
as well as some phase factors.

%The language of corepresentations and covariant algebras was used
%to formulate symmetry properties in the quantum field theory \cite{11,12},
%so that additional
%non-commuting objects from quantum groups did not enter into final
%physical quantities.

%The notion of universal $R$-matrix for the quasi-triangular Hopf algebras was
%given in \cite{1} as a canonical element combined of dual basis of dual Hopf
%algebras. This object is very useful, in particular for tensor products of
%Hopf algebra tensor operators \cite{13}. Corresponding universal analog of
%the quantum group matrix $T$ \cite{3} as a canonical
%element of $(F)^* \otimes F$ where $(F)^* $
%is the quantum algebra was analyzed in \cite{14}. Such
%universal objects are of importance to formulate quantum group invariant
%physical systems \cite{10,15} and/or quantum group Fourier
%transform \cite{16}.

2.5 The next example of a covariant system is related to the
reflection equation
algebra ${\cal K}$ (or the q-Minkowski space-time algebra, or the quantum
sphere algebra) (see e.g. \cite{16}). Its quantum group covariance depends
on the set of $R$-matrices in the defining equation (a reflection equation)

\begin{equation}
R_{12}^{(1)} K_1R_{12}^{(2)}K_2 = K_2R_{12}^{(3)}K_1R_{12}^{(4)},
\end{equation}

\noindent
with the coaction $\varphi (K) = K' = T K S $ where
$R^{(j)},\; j = 1,\; \ldots , 4$ define the commutation relations of the $T$
and $S$ entries \cite{16}. For the simple $SU_q(2)$ covariant case one has
$R_{12}^{(1)} = R_{12}^{(3)} =R_{12},\quad
R_{12}^{(2)} = R_{12}^{(4)} = R_{21}$ with
the $sl_q(2)$ matrix $R$ and $T=S^{\dagger}=S^{-1}$. $K$ is the following
$2 \times 2$ matrix of generators
$$
K =
\left(
\begin{array}{cc}
\alpha & \beta \\
\gamma & \delta
\end{array}
\right).
$$

\noindent
So the algebra ${\cal K}$ has four generators:
$\alpha,\; \beta,\; \gamma,\; \delta$ with relations
$$
\begin{array}{ll}
\alpha \beta = q^{-2} \beta \alpha, & [\delta, \beta] = q^{-1} \lambda
\alpha \beta, \\
\alpha \gamma = q^2 \gamma \alpha, & [\delta, \gamma] = -q^{-1}
\lambda \gamma \alpha, \\
{}[\alpha, \delta] = 0, & [\beta, \gamma] = q^{-1} \lambda
(\delta - \alpha) \alpha,
\end{array}
$$

\noindent
and two central elements
$$
c_1 = q^{-1} \alpha + q \delta, \quad
c_2 = \alpha \delta - q^2 \beta \gamma.
$$

One has the covariance of ${\cal K}$ with respect to the quantum group
$S U_q (2)$ with the coaction $\varphi: {\cal K} \rightarrow S U_q (2)
\otimes {\cal K}$ which is easy to write using the matrix form
$$
\varphi (K) = K' = U K U^{\dagger},
$$

\noindent
where $U = (U^{\dagger})^{-1}$ is the following
$2 \times 2$ matrix of the $S U_q (2)$ generators
$$
U =
\left(
\begin{array}{cc}
a & q b \\
-b^{\dagger} & a^{\dagger}
\end{array}
\right) .
$$

Due to the fact that the q-determinant is equal to one
$$
a a^{\dagger} + q^2 b^{\dagger} b = a^{\dagger} a + b^{\dagger} b = 1
$$

\noindent
the quantum group $S U_{q}(2)$ has essentially one unitary irreducible
representation ${\cal H}_{F}$ \cite{VS} with vacuum
state $|0>: a|0>= 0,\; b|0> = 1|0>,$
$$|n> = (a^{\dagger})^n |0> /\Pi_{j} c_j,  \quad c^2_n = (1 - q^{2n}). $$

The algebra ${\cal K}$ with the *-operation $K = K^{\dagger}$ has
many irreducible representations \cite{17}. Let us fix
one of them ${\cal H}_1$ then after the coaction $\varphi$ the transformed
algebra $K'$ generated by $(K')_{ij} = \varphi (K_{ij})$ is defined in the
tensor product ${\cal H}_F \otimes {\cal H}_1$. Hence, there is the problem of
the tensor product decomposition on the irreducible representations. The
transformed generators look rather cumbersome in terms
of the original generators
$$
K' =
\left(
\begin{array}{cc}
a & q b \\
-b^{\dagger} & a^{\dagger}
\end{array}
\right)
\left(
\begin{array}{cc}
\alpha & \gamma^{\dagger} \\
\gamma & \delta
\end{array}
\right)
\left(
\begin{array}{cc}
a^{\dagger} & -b \\
q b^{\dagger} & a
\end{array}
\right),
$$

\noindent
hence $e.g.$ $ \varphi (\alpha) = a a^{\dagger} \alpha + q b a^{\dagger}
\gamma + q a b^{\dagger} \gamma^{\dagger} +
q^2 b b^{\dagger} \delta.$

Let us consider the very simple (one-dimensional) irreducible
representation of the algebra
${\cal K}: \alpha = \delta = 0, \quad \gamma \in {\bf R}$.
The factor ${\cal H}_1$ is one-dimensional and ${\cal H}_F$
has to be decomposed into the irreducible representations
of $\cal K$. To reach this aim one has to find eigenvalues $\lambda$ and the
corresponding eigenvectors $|\lambda >$ of
$\alpha' = \varphi (\alpha) = q \gamma (b a^{\dagger} + a b^{\dagger})$
in ${\cal H}_F$.
Those of them related by $| \lambda_{n + 1} > \simeq \varphi (\gamma )
| \lambda_n >, \quad \lambda_{n + 1} = q^2 \lambda_n$
give rise to the invariant subspace of ${\cal H}_F \; w.r.t. \; {\cal K}'$.
The Hermitian operator $(b a^{\dagger} + a b^{\dagger})$ is a Jacobian matrix
with the entries $q^n c_n$ on the sub-diagonal. Hence the problem of the
non-trivial deficiency indices could take place \cite{14}.

2.6 A more complicated covariant system is related to the quantum super group
$OSp_q(1|2)$ \cite{2}. According to the general arguments of the Introduction,
the coaction map gives rise to the extension of the dynamical system and to the
representation of the covariantly transformed system in the tensor product
with one of the factor being an irreducible representation of the
corresponding QG. To find an irreducible unitary representation
of the quantum super-group $O S p_{q} (1|2)$
one has to introduce a *-operation and to analyse the commutation relations
among the generators $T_{ij}, i,j=1, 2, 3$. The matrix $T$ of
the $O S p_{q} (1|2)$ generators is even and has the dimension 3 in the
fundamental representation and the grading (0, 1, 0). The compact form of
the quadratic relations among the generators is given by the
$Z_2$-graded  FRT-relation ( $Z_2$-graded tensor product [2, 3])
$$
\hat RT \otimes T = T \otimes T \hat R \;.
$$

The $osp (1|2)$- $R$-matrix $\hat R$ has the spectral decomposition [2]
$$
\hat R = q P_5 - q^{-1}P_3 - q^{-2}P_1 \;,
$$

\noindent
where the projector indices refer to the dimension $4s+1$ of the subspaces
corresponding to spin 1, 1/2, 0 (see their explicit expressions in [9]).

Due to the structure of the $R$-matrix and the orthosymplectic condition
$T^{st}C_qT= \gamma C_q$ [2], there are only three independent
generators among
nine entries of $T$. One can easily see this from
the Gauss decomposition \cite{15} of the matrix $T$
$$
T = T_L T_D T_U\;,
$$

\noindent
where $T_L,\; T_U$ are lower and upper triangular matrices with unities on
their diagonal and the diagonal factor $T_D = diag (A,\; B,\; C)$. Among the
Gauss decomposition generators one finds three independent ones:
$A,\; (T_L)_{21}$ and $(T_U)_{12}$, while the element $B$ is central
\cite {15}.
Introducing the nine elements
of the Gauss decomposition: $T_D = diag (A, B, C), (T_L)_{21, 31, 32} =
(x, y, z)$ and $(T_U)_{12, 13, 23} = (u, v, w)$ one finds from the FRT-relation
\cite{2}:
$$ A = T_{11}, \quad x = T_{21} / A, \quad y = x^2 / \omega , \quad
z = x / q^{1/2},$$
$$
u = T_{12} / q A, \quad v = u^2 / \omega , \quad w = -q^{1/2} u, \quad
B=T_{22}- T_{21}(T_{11})^{-1} T_{12}\;,
$$

\noindent
where  the elements $B = AC = CA$ are central and
$\omega = q^{1/2} - q^{-1/2}$.

Due to the commutativity of $T_{13}$ and $T_{31}$ which are conjugated to each
other  $T_{31} = -T^{\dagger}_{13} / q$
according to the *-operation from \cite{2},
these elements are diagonal in the $Z_2$-graded Fock representation
${\cal H}_F$
with the vacuum: $T_{21} | 0 > = 0$ and the element $T_{12}$ as a creation
operator.

{}From the structure of the quadratic relations among
the generators $T_{ij}$ \cite{2} it follows that
the four elements  $T_{12}, \; T_{32}, \;
T_{13}, \; T_{31}$ form a subalgebra  of the $OSp_q(1|2)$
$$
T_{13} T_{12} = q^{-1} T_{12} T_{13} \; , \;
T_{13} T_{32} = q T_{32} T_{13} \; , \;
T_{13} T_{31} = T_{31} T_{13}\;,
$$
$$
T_{12} T_{32} + q T_{32} T_{12} = \lambda q^{1/2} T_{31} T_{13} \;.
$$

\noindent
Hence, the irreducible
representation in the Fock space ${\cal H}_F$ is given by $T_{12}$
as creation operator and
$T_{32}$ as annihilation operator, while
$$
T_{11}= q (T_{12})^2/ (\omega T_{13})
$$

\noindent
with $T_{13}$ and $T_{31}$ being diagonal in the basis
$|n> \simeq (T_{12})^n |0>$.

Let us now define the quantum $OSp$-plane, which is an associative super
algebra ${\cal A}$ with three generators $a,\; \xi ,\; b$ and
the $Z_2$-grading
$p(a) = p (b) = 0,\; p (\xi) = 1$.
Taking into account the similarity of the quantum
super-group $O S p_{q} (1|2)$ to the symplectic group case [1, 7] the defining
relations of ${\cal A}$ can
be written in the $R$-matrix form with a central extension

\begin{equation}
\hat R X \otimes X = q X \otimes X + c_2 J \;,
\end{equation}

\noindent
where $\hat R$ is the $osp_{q} (1|2)$ $R$-matrix, $X = (a, \xi, b)^t $
and the nine component vector
$$J = (0, 0, -q^{-1/2}, 0, 1, 0, q^{1/2}, 0, 0)$$
is the rewritten invariant matrix $C_q$ \cite{2}. One has for the generators
the quadratic relations

\begin{equation}
\xi a = a \xi /q \; , \; \xi b = q b \xi \;,\; [a, b] = \mu \xi^2 \;.
\end{equation}

\noindent
where $\mu = q^{1/2} + q^{-1/2} = \lambda / \omega $.
The vector $J$ is the eigenvector
of the rank one projector $P_1$ which gives rise to the centrality of the
element
$$
c_2 =\lambda (q\xi^2 /\omega - ab)/q^2 (q^{3/2}+ q^{-3/2}) =
 \lambda ( \xi^2 / q\omega - ba)/q^2 (q^{3/2}+ q^{-3/2}) .
$$
\noindent
This central element is invariant under the $O S p_{q} (1|2)$ coaction:
$X \rightarrow T X$. The algebra $\cal A$  was identified
in [9] as a twisted q-super-oscillator. Although $\cal A$ has the same
irreps as (3) with $b=a^{\dagger}$ the coaction is more complicated with
respect to (11) of the Subsec.2.2 for it includes the number or scaling
operator $ \xi = \eta q^N$ as well.

\section{Conclusion}

The problems of the quantum group coaction interpretation and the corresponding
tensor product decomposition are especially interesting in the framework of
the Poincare group deformation \cite{L,10,11}. The corresponding quantum
group has many generators and rather complicated quadratic relations
among them. Even in the very simple (trivial ?) case when the deformation
of the Poincare group is given by the twisting \cite{10} there are two Weyl
generators defining the representation. The Hamiltonian of the relativistic
system being only covariant under the group transformation law will get
extra degrees of freedom after the quantum group coaction \cite{11}. Another
kinematical group: the $q$-Galilei algebra $G_q$, was connected with
the $XXZ$-model
dispersion relation due to the equivalence of the trigonometric function
addition law and a non-commutative coproduct \cite{EC}. Realizing the
generators of $G_q$ in terms of the local spin operators one can obtain by
the duality the quantum group coaction.  Further interesting possibilities
for the representation theory refer to the case when coproduct or coaction
maps the original algebra into a tensor product with non-commutative factors
(see e.g. \cite{16,18}).

{\bf Acknowledgement}.
The authors thank P. Presnajder, R. Sasaki and \\
M. Scheunert for useful discussions.
PPK appreciate
the influence and support of the Non-perturbative quantum field theory
workshop of the Australian National University and Laboratoire de Physique
Theorique et Haute Energie associe au C.N.R.S.


\newpage

\begin{thebibliography}{99}
\bibitem{1}L. D. Faddeev, N. Yu. Reshetikhin and L. A. Takhtajan, Alg.
Analiz 1 (1989) 178 ( Leningrad Math. J. 1 (1990) 193 ).
\bibitem{2} P. P. Kulish and N. Yu. Reshetikhin, Lett. Math. Phys. 18 (1989)
143.
\bibitem{3} M. Chaichian and P. P. Kulish, Phys. Lett. B233 (1990), 72.
\bibitem{4} M. Chaichian, P. P. Kulish and J. Lukierski, Phys. Lett.
B262 (1991) 43.
\bibitem{5} W. Pusz and S. L. Woronowicz, Rep. Math. Phys. 27 (1989) 231.
\bibitem{6} V. Rittenberg and M. Scheunert, J. Math. Phys. 33 (1992) 436.
\bibitem{7} P. P. Kulish, Phys. Lett. A161 (1991) 50.
\bibitem{8} P. P. Kulish, Theor. Math. Phys. 85 ( 1991) 157;
i.b. 94 (1993) 193.
%Theor. Math. Phys.
\bibitem{9} F. Thuillier and J. C. Wallet, Phys. Lett. B323 (1994) 153.
\bibitem{V} P. P. Kulish, Proc. of Varenna School on Quantum Groups:
Theory and Applications, to be published by Plenum Press, 1995.
\bibitem{L} J. Lukierski, A. Nowicki, H. Ruegg and V. Tolstoy,
Phys. Lett. B264 (1991) 331.
\bibitem{10} M. Chaichian and A. Demichev, J. Math. Phys. 36 (1995) 398.
\bibitem{11} J. A. de Azc\'arraga, P. P. Kulish and F. R\'odenas, Phys. Lett.
B351 (1995) 123; hep-th/9405161, Fortschr. Phys. (1995).
\bibitem{12} M. Arik, In: Symmetries in science VI, Plenum Press, 1993, 47.
\bibitem{13} M. Chaichian, H. Grosse and P. Presnajder,
J. Phys. A27 (1994) 2045.
\bibitem{14} I. M. Burban and A. U. Klimyk, Lett. Math. Phys. 29 (1993) 13.
\bibitem{15} E. V. Damaskinsky, P. P. Kulish and M. A. Sokolov, Zap. Nauch.
Semin. POMI, 224 (1995) 155; ESI-95-217; q-alg/ 9505001.
\bibitem{16} P. P. Kulish and R. Sasaki, Prog. Theor. Phys. 89 (1993) 741.
\bibitem{17} P. P. Kulish, Alg. Analiz, 6 (1994) 195.
\bibitem{EC} F. Bonechi, E. Celeghini, R. Giachetti, E. Sorace and
M. Tarlini, J. Phys. A 25 (1992) L939; Phys. Rev. B 46 (1992) 5727.
\bibitem{18} T. H. Koornwinder, In: Orthogonal polynomials:
theory and practice, (ed. P.Nevai) NATO AS Serie 1991, 257.
\bibitem{B} O. Babelon and D. Bernard, Commun. Math. Phys. 149 (1992) 279.
\bibitem{VS} L. Vaksman and Ya. Soibelman, Func. Analiz Pril. 22:3 (1988) 1;\\
Ya. Soibelman, Int. J. Mod. Phys. 7 Supp. 1B (1992) 859.

\end{thebibliography}
\end{document}